\documentclass[runningheads]{llncs}
\usepackage[T1]{fontenc}
\usepackage{amsfonts}
\usepackage{amssymb}
\usepackage{graphicx}
\begin{document}
\title{Uncertainty, Imprecise Probabilities and Interval Capacity Measures on a Product Space}
\titlerunning{Uncertainty and Imprecise Probabilities}
\author{Marcello Basili\inst{1} \and
Luca Pratelli\inst{2}}
\authorrunning{Basili and Pratelli}
%
\institute{DEPS University of Siena, Italy,\ \email{marcello.basili@unisi.it}\\ \and Naval Academy, Leghorn, Italy,\ \email{luca\_pratelli@marina.difesa.it}}
\maketitle
\begin{abstract}
In Basili and Pratelli (2024), a novel and coherent concept of interval probability measures has been introduced, providing a method for representing imprecise probabilities and uncertainty. Within the framework of set algebra, we introduced the concepts of weak complementation and interval probability measures associated with a family of random variables, which effectively capture the inherent uncertainty in any event.

This paper conducts a comprehensive analysis of these concepts within a specific probability space. Additionally, we elaborate on an updating rule for events, integrating essential concepts of statistical independence, dependence, and stochastic dominance.
\end{abstract}
\section{Introduction}

The Theory of Belief Functions, initially introduced by Dempster (1967) to extend Bayesian statistical inference, was further developed by Shaffer (1976) to tackle uncertainty representation. This evolution led to Uncertain Evidence Theory, where the combination rule became pivotal. Walley (1991) emphasized the necessity of coherent lower and upper probability measures for a comprehensive theory of imprecise probabilities. 

Belief functions enable the representation of incomplete knowledge, the pooling of evidence, and the updating of beliefs in the face of new evidence. Denoeux (2019) outlined decision-making criteria within the belief functions framework. Dubois and Denoeux (2012) delineated two methods for conditioning belief functions: the Dempster rule for revising a plausibility function and a method tailored for prediction based on observation. Additionally, Dubois et al. (2023) demonstrated that a straightforward belief function logic can be derived by integrating Lukasiewicz logic into minimal epistemic logic, forming a probabilistic logic within a modal logic framework

This concise paper revisits the concept of weak complementation, a distinctive notion of set-theoretic complementation within the algebra of sets (or events) of a generic non-empty set, as discussed in Basili and Pratelli (2024). This framework offers a natural and coherent approach to conceptualizing uncertainty linked with an event. By interpreting eventualities as the causal factors leading to the occurrence of a given random phenomenon, we introduce the set of indecisive eventualities (uncertain opportunities) and quantify the degree of uncertainty associated with an event within a specific probability space of particular interest in applications. Furthermore, within this framework, we elaborate on fundamental concepts such as imprecise probabilities, interval capacity measures and conditional interval probability measures. Finally, we illustrate the notion of stochastic dominance between two random variables based on interval probability measures and exhibit an example of the product of interval probability measures.

\section{ A Particular Probability Space}
Consider the product space $\Omega$ defined as $E\times\{0,1\}^n$, where $E$ is a countable space and $n$ is a natural number greater than 0. In this space, an element $(x,\omega_{i_1\ldots i_n})$ is an eventuality, where $x$ belongs to $E$, and $\omega_{i_1\ldots i_n}$ denotes a sequence of outcomes or causal factors, each represented by 0 or 1. Subsets $H$ of these eventualities are events and contribute to the determination of a given random phenomenon. However, due to imprecision in the elements of $\{0,1\}^n$, indecisive eventualities arise when considering a specific event $H$.  For instance, when $n=2$ and $E$ is a singleton $\{x_0\}$, the famous example of the umbrella (Keynes 1921, 28) has been studied on $\{0,1\}^2$ (equivalent to $E\times\{0,1\}^2$) where $\{\omega_{00},\omega_{11}\}$ is the set of indecisive eventualities of $H=\{\omega_{10}\}$ since the negation of a given cause does not necessarily imply the non-occurrence of the represented event. (See Example 5 of Basili and Pratelli (2024))

\subsection{Incompatible Eventualities}
In the following scenario, any pair of distinct eventualities $$(x,\omega_{i_1\ldots i_n}),(x,\omega_{j_1\ldots j_n})$$ where $0<|i_1-j_1|+\ldots+|i_n-j_n|<n$ is considered not incompatible due to the imprecision or indecision represented by the elements of $\{0,1\}^n$. However, denoted by $\omega_{i_1^*\ldots i_n^*}$ as the (unique) negation of $\omega_{i_1\ldots i_n}$ if $|i_1-i_1^*|+\ldots+|i_n-i_n^*|=n$, the two eventualities $$(x,\omega_{i_1\ldots i_n}),(y,\omega_{i_1^*\ldots i_n^*})$$ can be considered incompatible for any pair $x,y$ of elements of $E$. It's observed that there are $2^{n-1}$ events of the type $\{\omega_{i_1\ldots i_n},\omega_{i_1^*\ldots i_n^*}\}$ and we denote by  $Z_j$ the events of the form $E\times \{\omega_{i_1\ldots i_n},\omega_{i_1^*\ldots i_n^*}\}$, with $j=1,\ldots,2^{n-1}$. Obviously, events $Z_j$ constitute a partition of $\Omega$. In the particular case of $n=1$, two distinct eventualities are always incompatible.

\subsection{Set of Indecisive Eventualities and Weak Complementation}

    To determine the set of uncertain opportunities or indecisive eventualities of an event $H$ based on the previously stated and elaborated concepts in Basili and Pratelli (2024), we define $$H_{ind}=\cup_{j:H\cap Z_j=\emptyset}\  Z_j$$ 	
 In simpler terms, $H_{ind}$ is the set of eventualities that are not incompatible with the elements of $H$ while $H^c\setminus H_{ind}$ is the set of incompatible elements with the eventualities of $H$. The set $H^c_w=H^c\setminus H_{ind}$ can be interpreted as the weak complement of $H$, as it exclusively comprises eventualities that unequivocally imply the non-occurrence of $H$. Thus, for any event $H$, it yields $$\Omega=H\cup H^c_w\cup H_{ind}.$$It is important to note that $K_{ind}\subset H_{ind}$ implies $K_{ind}\subset H_{ind}$ but it does not imply $K_{w}^c\subset H_{w}^c.$
 
\subsection{Degree of Uncertainty Associated with an Event}
For any event $H$,  a degree of uncertainty can be associated with any eventuality incompatible with elements of $H$. Specifically, considering a generic mapping $r:\Omega\mapsto [0,1]$ representing the degree of uncertainty on $\Omega$, the random variable $Y_H$ defined by $$Y_H(x,\omega_{i_1\ldots i_n})=I_{H_{ind}}(x,\omega_{i_1\ldots i_n})r(x,\omega_{i_1\ldots i_n})$$ is the $r$-uncertainty associated with $H$. It is agreed that $Y_H=0$ when $H_{ind}=\emptyset$. In the particular case of $r=1$, it follows $Y_H=I_{H_{ind}}$ whereas if $r=0$ the $r$-uncertainty of $H$ is always null. In certain aspects, $Y_H(x,\omega_{i_1\ldots i_n})$ is akin to the degree of indeterminacy of $(x,\omega_{i_1\ldots i_n})$ concerning $H$ in intuitionistic fuzzy sets (see Atanassov, 1986).

\section{Imprecise Probabilities}
 Now, let's consider a probability measure $P$ defined on the events of $\Omega$ and a degree of uncertainty $r$ defined on $\Omega.$ The interval probability measure $Q_r$ defined by
 $$Q_r(H)=\big[P(H),P(H)+E[rI_{H_{ind}}]\big]$$ is said to be the {\it imprecise probability measure associated with $P$ and $r$}. When $H_{ind}$ is not $P$-negligible, $Q_r(H)$ coincides with $[P(H),P(H)+P(H_{ind})E_{P_{H_{ind}}}[r]]$. It's worth recalling (see [2]) that a generic interval probability measure (or imprecise probability) $Q$ is a function from the events of $\Omega$ to the closed subintervals of $[0,1]$, satisfying:
 
\begin{itemize}
\item The left-end extreme of $Q$ is a probability measure.
\item For every pair of events $H_1,H_2$, with $H_1\subseteq H_2$, it holds $|Q(H_2)|\leq |Q(H_1)|$, where $|H_i|$ denotes the width of $H_i$.
\end{itemize}

\smallskip\noindent When $r=1$, for any event $H$ it holds $$Q_1(H)=[P(H),P(H)+P(H_{ind})].$$ Notwithstanding, since $P(H) + P(H_{ind}) = 1 - P(H^c_w)$, the right-end extreme of $Q_1(H)$ consequently aligns with the left-end extreme of $Q_1((H^c_w)^c)$, where  $(\cdot)^c$ denotes the usual set complementation. Finally
$$|Q_r(H)|\leq P(H_{ind})=\sum_{\omega_{i_1\ldots i_n},\omega_{i_1^*\ldots i_n^*}\in H^c}\big(f(\omega_{i_1\ldots i_n})+f(\omega_{i_1^*\ldots i_n^*})\big)$$ where $f$ is the probability mass function defined by $f(\omega_{j_1\ldots j_n})=P(E\times\{\omega_{j_1\ldots j_n}\})$.
 
 \subsection{Interval Capacity Measures}
Instead of the probability measure $P$, a more generalized definition of an interval measure can be derived by considering a capacity $\nu$, which is a function defined on the events of $\Omega$ satisfying  $\nu(\emptyset)=0, \nu(\Omega)=1$ and $\nu(H)\leq\nu(K)$ when $H\subseteq K$. More
precisely, $Q_r$ could be defined as follows 
$$
Q_r(H)=\big[\nu(H),\nu(H)+\int_0^1 \nu(H_{ind}\cap \{r\geq t\})\, dt\, \big] 
\cap[0,1]$$ where $\int_0^1 \nu(H_{ind}\cap \{r\geq t\})\, dt$ represents the Choquet integral of $rI_{H_{ind}}$. Moreover, if $\nu$ is a super-additive set function, it holds that:
$$Q_r(H)\subseteq  \big[\nu(H),\int_0^1 \nu(H\cup (H_{ind}\cap \{r\geq t\}))\, dt\, \big]$$ and $|Q_r(K)|\leq |Q_r(H)|$ when $H\subseteq K$. In this context, the interval capacity measure of $H$ could be considered $\big[\nu(H),\int_0^1 \nu(H\cup (H_{ind}\cap \{r\geq t\}))\, dt\, \big]$. However,  the map $$Q^\prime_r:H\mapsto \big[\nu(H),\int_0^1 \nu(H\cup (H_{ind}\cap \{r\geq t\}))\, dt\, \big]$$ does not necessarily satisfy $H\subseteq K\Longrightarrow |Q^\prime_r(K)|\leq |Q^\prime_r(H)|$ unlike $Q_r$.

\subsection{Conditional Interval Probability Measures}

Let $H$ be a non-negligible event ($P(H)>0$). An interval probability measure can be defined to express the degree of confidence in realizing the random phenomenon conditioned on the causes that determined $H$. For each event $A$, we define the \textit{probability of $A$ conditioned on $H$ relative to $Q_r$} as
$$Q_r(A|H)=\big[{\frac{{P(A\cap H)+E[rI_{A\cap H_{ind}}]}}{{P(H)+E[rI_{H_{ind}}]}}},{\frac{{E^{P}[(I_A+rI_{A_{ind}})(I_H+rI_{H_{ind}})]}}{{P(H)+E[rI_{H_{ind}}]}}}\big].$$ $Q_r(\cdot|H)$ is an imprecise probability in the sense of $[2].$ If $r=1$, it leads to
$$
Q_1(A|H)=\Big[P\big(A\, | (H_w^c)^c),P\big((A_w^c)^c\, |(H_w^c)^c\big)\Big] 
$$It is noteworthy that $Q_1$ possesses left and right extremes that adhere to the duality rule concerning weak complementation. Moreover, it holds
$$ Q_1(\Omega|H)=Q_{1}((H_w^c)^c\, |H)=[1,1],\qquad Q_1(H\, |H)=[\frac{P(H)}{1-P(H_w^c)},1].$$

\smallskip\begin{remark}
Let $\nu$ be a super additive capacity  i.e. for any pair of disjoint events $A,B$, $\nu$ satisfies $\nu(A)+\nu(B)\leq \nu(A\cup B)$.  If $\nu(H^c)\neq 1$, Dempster and Shafer define  $$\nu(A|H)={{\nu((A\cap H)\cup H^c)-\nu(H^c)}
\over{1-\nu(H^c)}}$$ for any event $A$. 
This definition does not equate to standard Bayesian updating if $\nu$ is not additive. However, if $H^c$ is replaced with the weak complementation $H^c_w$, and thus $H$ replaced with $(H^c_w)^c$, interpreting the plausibility $1 - \nu(H^c_w)$ as the measure of $(H^c_w)^c$ according to a probability $P$, we can see that $${{\nu((A\cap (H^c_w)^c)\cup H^c_w)-\nu(H^c_w)}
\over{1-\nu(H^c_w)}}={{\nu((A\cup H^c_w)-\nu(H^c_w)}
\over{1-\nu(H^c_w)}}$$ can be understood as the ratio of measures of $A \cap (H^c_w)^c$ and $(H^c_w)^c$ according to $P$. In other words, the left endpoint of $Q_1(A|H)$. Similarly, the interpretation of the right endpoint of $Q_1(A|H)$ holds true. Therefore, $Q_1(A|H)$ can represent the conditional interval capacity measure associated with $\nu$ since the extremes are  two similar conservative belief and plausibility degrees conditional to $H$ by Dempster and Shafer.
\end{remark}

\smallskip\begin{remark} Is there a 'reasonable' way to define the concept of conditional probability with respect to $H$ according to a general interval capacity measure? In Basili and Pratelli (2024), we proposed only a hypothetical and partial definition, deferring the concept to further work (See also Denoeux et al (2020)). Within this framework, inspired by the notions introduced in section 3.1, if $\nu$ is a super-additive capacity and $H$ is an event with $\nu(H)>0$, the definition of $Q_r(A|H)$ can be expressed as  
$$Q_r(A|H)=\big[ \frac{I(A)}{I(\Omega)},\frac{I(A)+J(A_{ind})}{I(\Omega)}\big]$$
where $I(B)= \int_0^1 \nu\big(B\cap (H\cup (H_{ind}\cap \{r\geq t\}))\big)\, dt$ and $$J(B)=\int_0^1 \nu\big(B\cap (H\cup H_{ind})\cap \{r\geq t\}\big)\, dt.  $$It is noteworthy that $$Q_r(A|H)\subseteq \big[ \frac{I(A)}{I(\Omega)},\frac{I(A)+I(A_{ind})}{I(\Omega)}\big]$$ and $A\subseteq B\Longrightarrow |Q_r(B|H)|\leq |Q_r(A|H)|$. By arguing as in section 3.1, one might introduce the  definition $$Q^\prime_r(A|H)=\big[ \frac{I(A)}{I(\Omega)},\frac{I((A^c_w)^c)}{I(\Omega)}\big]$$ noting that $$Q_r(A|H)\subseteq \big[ \frac{I(A)}{I(\Omega)},\frac{I(A)+I(A_{ind})}{I(\Omega)}\big]\subseteq Q^\prime_r(A|H).$$ Unfortunately, $Q^\prime_r(\cdot|H)$, like $Q^\prime_r(\cdot),$ does not necessarily satisfy $A\subseteq B\Longrightarrow |Q^\prime_r(B|H)|\leq |Q^\prime_r(A|H)|.$ 
\end{remark}

\section{Interval Distribution and Stochastic Dominance}

As in Basili and Pratelli (2024), we consider a real-valued random variable $X$ and define the {\it
distribution function (according to a generic imprecise probability $Q_r$) of $X$} as the function $F$ from the real line to the set of closed
subintervals of $[0,1]$ given by 
$$
F(t)=Q_r(X\leq t)
$$The widths of the distribution function $F$ depend significantly on the values of $Z_j$, as demonstrated in the following example.

\begin{example}
A special case arises when $X(Z_j)=t_j$. In this case, it holds that 
\[
F_{r}(t)=\Big[P(X\leq t),P(X\leq t)+(1-P(X\leq t))E_{P_{\{X>t\}}}[r]\Big]
\]
where $P(X\leq t)=\sum_{j:t_j\leq t} P(Z_j)$ and $E_{P_{\{X>t\}}}[r]$ is the mean of $r$ with respect to $P(\cdot\,|\,X>t)$. In particular, $$F_1(t)=\Big[P(X\leq t),1\Big].$$
 If $t_j$ is non-decreasing and $
Y(Z_j)$ is a countable subset of $]t_{j-1},t_j]$, then $X\geq Y$ and the distribution function $G_1$,
according to $Q_{1}$ is
given by 
\[
G_1(t)=\Big[P(Y\leq t), 1-P(Z_{i_*})\delta_ t\Big], 
\]
where $t_{i_*-1}\leq t<t_{i_*}$, and  $\delta_t=1$ if $\{t_{i_*-1}<Y< t\neq\emptyset\}$ and $\delta_t=0$ elsewhere. It's noteworthy that the
function representing the right endpoint of $G_{1}$ is neither
increasing nor decreasing. Additionally, it holds $P(X\leq t)\leq P(Y\leq t)$ and 
$|F_1(t)|\geq |G_1(t)|$ for any $t\in {\mathbb{R}}$.
 
\end{example}
 
 \smallskip These considerations suggest a notion of interval stochastic dominance. More precisely, if $F$ and $G$ are the two (interval) distribution of $X$ and $Y$, it is said that
{\it $X$ stochastically dominates $Y$} if 
$$F_l(t)\leq G_l(t) \ \ \ {and}\ \ \ |G(t)|\leq
|F(t)|\quad \ \ \ \ \forall t\in {\mathbb{R}}$$
where $F_l$ and $G_l$ are the left endpoints of $F$ and $G$.

 \smallskip\noindent A similar definition can be extended when $X$ and $Y$ are random vectors. This concept of interval stochastic dominance enables the comparison of random variables in terms of uncertainty and provides a straightforward guideline for making choices among various alternatives. Finally, these results can be generalized when $Q$ is an interval capacity measure obtained by a super additive capacity $\nu$. However, if $\nu$ is only a sub-additive capacity, the stochastic dominance cannot be extended because $|G(t)|\leq |F(t)|$ is not necessarily true when $X\geq Y$.

\section{Product of interval probabilities measures}

Let's consider the space $E\times\{0,1\}^n$ for simplicity when $n=2$, and $E$ is a singleton $\{x_0\}$.  On this space,  let $Q_1$ be the imprecise probability measure  associated with $P$ (with $r=1$). It's worth noting that the product space $\Omega=(E\times\{0,1\}^2)^2$ is obviously equivalent to $E\times\{0,1\}^4$. We define the interval product measure $Q_1\otimes Q_1$ on the events $H$ of $\Omega$ as
$$Q_1\otimes Q_1(H)=[P\otimes P (H),P\otimes P(H)+P\otimes P(H^\prime_{ind})] $$ where $P\otimes P$ is the usual product measure and $H^\prime_{ind}$ is the subset of the eventualities of $\Omega$ that are not incompatible with the elements of $H$ with respect to the partition $$W_1=Z_1\times Z_1, \quad W_2=Z_2\times Z_2,\quad W_3=Z_1\times Z_2,\quad W_4=Z_2\times Z_1$$ where $Z_1=\{\omega_{00},\omega_{11}\}$ and 
$Z_2=\{\omega_{10},\omega_{01}\}$. Specifically, $$H^\prime_{ind}=\cup_{i:H\cap W_i=\emptyset}\  W_i$$
Observe that $Q_1\otimes Q_1$ does not coincide with any measure $Q^\prime_r$ associated with $P\otimes P$ and $r$ on $\Omega$ because, for $Q^\prime_r$, the event $H_{ind}$ is determined by eights disjoint events. In particular, if $H=H_1\times H_2$ with $H_1,H_2\subseteq E\times\{0,1\}$
then $$Q_1\otimes Q_1(H_1\times H_2)=[P(H_1)P(H_2),P(H_1)P(H_2)+P\otimes P(H^\prime_{ind})] $$
Since $P\otimes P(H^\prime_{ind})\neq P((H_1)_{ind})P((H_2)_{ind})$, the product measure $Q\otimes Q$ cannot be obtained by simply multiplying the endpoints and the widhts of intervals $Q_1(H_1)$ and $Q_1(H_2)$. It can be demonstrated that $$Q_1\otimes Q_1(H_1\times H_2)\subseteq Q^\prime_1(H_1\times H_2).$$ For example, if $H=\{(w_{10},w_{10})\}$, we have
$$Q_1\otimes Q_1(H)= [P(\{(w_{10}\})^2,P(\{(w_{10}\})^2+1-P(Z_2)^2]$$ which does not coincide with $ Q^\prime_1=[P(\{(w_{10}\})^2,1-P(\{w_{01}\})^2]$ because $$P(\{(w_{10}\})^2+1-P(Z_2)^2=1-P(\{w_{01}\})^2-2P(\{w_{01}\})P(\{w_{10}\}).$$ These results could be useful for studying the example of the umbrella with two observers. The two distinct measure $Q_1\otimes Q_1$ and $Q_1^\prime$ are two interval probability measures derived from the same $P\otimes P$ but through two different ways of evaluating uncertainty. In conclusion, it can be easily noted that there is a number of ways of evaluating uncertainty which grows exponentially with the dimension $n$.

\end{document}